\documentclass[12pt]{article}
\usepackage{amsfonts,amssymb,version}

\def\AA{{\bf A}}

\def\bb{{\bf b}} 
\def\C{\mathbb C}

\def\D{{\mathcal D}}
\def\E{{\mathcal E}}
\def\EE{{\bf E}} 
\def\F{{\mathcal F}}
\def\FF{\mathbb F}

\def\GG{{\mathbb G}}

\def\L{{\mathcal L}}
\def\LL{{\bf L}} 
\def\M{{\mathcal M}}
\def\N{{\mathcal N}}

\def\OO{{\mathcal O}}

\def\Qq{{\mathcal Q}}

\def\R{\mathbb R}
\def\T{{\cal T}}
\def\XX{\mathfrak X}
\def\YY{\mathfrak Y}
\def\Z{\mathbb Z}
\def\Zz{{\mathcal Z}}

\def\coker{\mathop{\rm coker}\nolimits}

\def\unhom{\underline{Hom}}
\def\id{\mathop{\rm id}\nolimits}

\def\ker{\mathop{\rm ker}\nolimits}

\def\Proj{\mathop{\rm Proj}\nolimits}

\def\qed{\hfill$\square$\medskip}
\def\rank{\mathop{\rm rank}\nolimits}

\def\spec{\mathop{\rm Spec}\nolimits}

\def\Sym{\mathop{\rm Sym}\nolimits}

\mathchardef\mhyphen="2D

\let\da\downarrow

\let\hra\hookrightarrow

\let\ov\overline

\let\un\underline
\let\upa\uparrow
\let\wh\widehat
\let\pa\partial

\newtheorem{theorem}{Theorem}
\newtheorem{proposition}[theorem]{Proposition}
\newtheorem{lemma}[theorem]{Lemma}

\def\stm{\refstepcounter{theorem}\paragraph{\thetheorem}}
\def\rem{\refstepcounter{theorem}\paragraph{Remark \thetheorem}}

\def\proof{\paragraph{Proof}}

\makeatletter \def\l@section{\@dottedtocline{1}{0em}{1.2em}} \makeatother

\begin{document}

\centerline{\Large\bf Quadric invariants and degeneration}
\centerline{\Large\bf in smooth-\'etale cohomology}

\bigskip

\centerline{\bf Saurav Bhaumik and Nitin Nitsure}

\bigskip

\centerline{\it Dedicated to M.S. Narasimhan on his completing 80 years}

\begin{abstract}
For a regular pair $(X,Y)$ of schemes over $\Z[1/2]$ 
of pure codimension $1$, we consider quadric bundles on $X$
which are nondegenerate on $X-Y$, but are minimally degenerate on $Y$. 
We give a formula for the behaviour of the 
cohomological invariants (characteristic classes) 
of the nondegenerate quadric bundle on 
$X-Y$ under the Gysin boundary map to the \'etale cohomology of $Y$ with
mod $2$ coefficients. 

The results here are the algebro-geometric analogs of 
topological results for complex bundles 
proved earlier by Holla and Nitsure, 
continuing further the algebraization program which was commenced
with a recent paper by Bhaumik. 
We use algebraic stacks and their smooth-\'etale cohomologies,
$A^1$-homotopies and Gabber's absolute purity theorem  
as algebraic replacements for the topological methods used earlier, 
such as CW complexes, real homotopies, Riemannian metrics and 
tubular neighbourhoods. Our results also hold for 
quadric bundles over algebraic stacks over $\Z[1/2]$.
\end{abstract}

\centerline{2010 Math. Sub. Class. : 14F20, 14D06, 14D23.}

{\footnotesize 

\tableofcontents

}

\pagestyle{empty}

\section{Introduction}

This paper extends to the algebraic category the topological results of 
Holla and Nitsure [H-N-1] and [H-N-2]. 
We address the following three questions,
which were answered in the earlier papers in
the topological category. 

Let $S$ be a scheme on which $2$ is invertible, 
for example, $S = \spec \Z[1/2]$. 
Let a nondegenerate quadratic triple $T= (E,L,b)$ of rank $n$ on a scheme 
$X$ over $S$ consist of a rank $n$ vector bundle $E$ and  
a line bundle $L$ on $X$, and 
a nondegenerate symmetric bilinear form $b$ on $E$ taking values in $L$.
As we explain later, 
such a triple $T$ is the same as a principal
$GO_{n,S}$-bundle where $GO_{n,S}$ is the general orthogonal 
group-scheme over $S$.
Let $BGO_{n,S}$ be its classifying algebraic stack,
and let $H^*(BGO_{n,S},\FF_2)$ denote its cohomology ring in the 
smooth-\'etale topology.
A characteristic class of degree $i$ for nondegenerate quadratic triples
of rank $n$ is an element $\alpha \in H^i(BGO_{n,S},\FF_2)$. Such an $\alpha$  
associates to any such $T$ an \'etale cohomology class 
$\alpha(T)\in H^i(X,\FF_2)$, which is functorial in the sense that 
it commutes with pull-backs under any $X'\to X$.


\pagestyle{myheadings}
\markright{Saurav Bhaumik and Nitin Nitsure: 
Quadric invariants and degeneration}


\noindent{\bf Question 1} What is the ring $H^*(BGO_{n,S},\FF_2)$ 
of characteristic classes for
nondegenerate quadratic triples of rank $n$ on schemes over $S$?

A nondegenerate quadric bundle $Q$ on $X$ of rank $n$ will mean an 
equivalence class of rank $n$ nondegenerate quadratic triples on $X$,
where we regard a triple $T= (E,L,b)$ as being equivalent to 
another triple $T\otimes K = (E\otimes K, L\otimes K^2, b\otimes 1_K)$,
where $K$ is any line bundle on $X$
(note that the vanishing of the quadratic forms 
associated to $b$ and $b\otimes 1_K$ 
define the same closed subscheme $C \subset P(E)=P(E\otimes K)$, 
which a bundle of nondegenerate quadrics).
A cohomological invariant of degree $i$
associates to any such $Q$ an \'etale cohomology class 
$\alpha(Q)\in H^i(X,\FF_2)$, which is functorial under pull-backs.
These form a subring $PH^*(BGO_{n,S},\FF_2) \subset H^*(BGO_{n,S},\FF_2)$.

\noindent{\bf Question 2} What is the ring $PH^*(BGO_{n,S},\FF_2)$
of cohomological invariants for
nondegenerate quadric bundles of rank $n$ on schemes over $S$?

We will now come to the third question, answering which is the main 
purpose of this paper. 
Let $X$ be a regular scheme on which $2$ is invertible, 
and let $Y\subset X$ be a closed regular subscheme of pure 
codimension $1$. 
We are interested in quadric bundles over $X$ which are nondegenerate
over $X-Y$ and minimally degenerate over $Y$. This means 
that if $Q$ is represented by a quadratic triple $T = (E,L,b)$,
then $b$ is nondegenerate of rank $n$ on $X-Y$, while rank of $b|_Y$ is 
$n-1$. 

\noindent{\bf Question 3} How do the cohomological invariants of 
the nondegenerate quadric bundle $Q|_{X-Y}$ behave under degeneration?

Before we come to describing our solution to the question 3, we make
some comments on questions 1 and 2. Let $QuasiProj_S$ be the 
category of all quasi-projective schemes over $S$. 
For any group scheme $G$ over $S$, we have a functor
$|BG|$ on $QuasiProj_S$, which associates to any $X$ the set of 
all isomorphism classes of \'etale locally trivial 
$G$-torsors on $X$. Let $H^*(-,\FF_2)$
be the functor on $QuasiProj_{S}$ which associates
to any $X$ the \'etale cohomology $H^*(X,\FF_2)$.
We show in section 4 (Proposition \ref{functorial char classes}) 
that when the base $S$ is an excellent Dedekind domain, 
for any reductive group-scheme $G$ over $S$, any 
cohomology class $\alpha \in H^*(BG,\FF_2)$ is the same as
a morphism of functors $\alpha'$ from $|BG|$ to $H^*(-,\FF_2)$
on the category $QuasiProj_{S}$.
Bhaumik [Bh] has recently shown this when the base $S$ is of the form 
$S= \spec k$ where $k$ is a field of 
characteristic $\ne 2$. All that remains for us to do is to point out how
the arguments in [Bh] -- which were based on an idea of Totaro [T] 
for approximating classifying stacks $BG$ using 
Geometric Invariant Theory and Jouanolou's trick -- 
extend from base $k$ to base $S$, using the 
extension of GIT to more general bases due to Seshadri [Se].

Our second comment on questions (1) and (2) is as follows.
The singular cohomology rings $R_n$ of the CW complexes 
$BGO_n(\C)$ for the complex Lie groups $GO_n(\C)$ with $\FF_2$ coefficients
have been explicitly 
determined in Holla and Nitsure [H-N-1] and  [H-N-2], in terms of 
generators and relations. The primitive subrings $PR_n\subset R_n$
have been determined in [H-N-2] in a closed form 
for all odd $n$, for $n=4$, and for all $n$ of the form $4m+2$. 
As shown in [Bh], for any separably closed base field $k$ of 
characteristic $2$,
we get the same answers:  $H^*(BGO_{n,k},\FF_2) = R_n$ and 
$PH^*(BGO_{n,k},\FF_2) = PR_n$. The situation over a general base,
such as $S = \spec \Z[1/2]$, is more complicated 
(see Remark \ref{when kummer simplifies}),
and will not be described here.

The answer to the question (3) is the main theorem of this paper 
(Theorem \ref{main theorem}), which 
describes how the \'etale cohomological invariants of the 
nondegenerate quadric bundle
over $X-Y$ behave under degeneration. Concretely, we describe how the 
quadric invariants behave under the Gysin boundary map 
$\pa : H^*(X-Y, \FF_2) \to H^{* -1}(Y,\FF_2)$   
from the \'etale cohomology of $X-Y$ to
the \'etale cohomology of the degeneration locus $Y$. 

An analogous result has been proved in [H-N-2]
over complex numbers  
in the topological category using singular cohomology. 
Some of the crucial techniques used in the topological case
-- namely, tubular neighbourhoods, orthogonal projections using 
Riemannian metrics and usual homotopy -- 
do not extend to the algebraic category. 
We replace the older arguments 
with the use of algebraic stacks which leads to a
much more natural argument -- one that  
also sheds light on the earlier topological result.
A by-product of the proof is that the main theorem works also for 
degenerating quadric bundles defined over 
regular pairs of algebraic stacks over $S$, 
so we have stated and proved it in this generality.

Extensive concrete calculations for the rings $R_n$ and $PR_n$
have been given in [H-N-1] and [H-N-2], and the images of the generators 
of $PR_n$ under  
the Gysin boundary map have been listed in all low ranks (up to rank 6) 
in [H-N-2]. As explained in section 5.5, these calculations 
now extend to the algebraic category
over separably closed fields of characteristic $\ne 2$.  
We expect this to be of some interest, 
given the important role that quadratic forms play in arithmetic and 
algebraic geometry. 

This paper is arranged as follows. In section 2, we recall basic facts about 
absolute purity over algebraic stacks, and prove a `multiplicity 
lemma' (Lemma \ref{multiplicity lemma}) 
for Gysin boundary maps which is useful later.
The section 3 introduces the basic notions about 
quadratic triples, quadric bundles and their degenerations, following 
[H-N-1] and [H-N-2]. 
We introduce the \'etale cohomological invariants of nondegenerate
quadratic triples and quadric bundles in section 4, and explain why
the arguments in [Bh] extend from base a field $k$ 
to an excellent Dedekind domain $S$,  
so as to show that the functorial characteristic
classes and cohomological invariants on quasi-projective schemes over such a 
base $S$ exactly form the rings $H^*(BGO_{n,S},\FF_2)$
and $PH^*(BGO_{n,S},\FF_2)$. 
The section 5 is devoted to the proof of our main theorem 
(Theorem \ref{main theorem}). The idea of the proof is to introduce 
a universal stack $\M$ for at most mildly degenerate quadratic triples 
(section 5.2) and to 
simulate a tubular neighbourhood of the degeneration locus $\D$ in $\M$,
up to $\AA^1$-homotopy, by another algebraic stack (introduced 
in section 5.3) which is the stackification of the `standard model' for 
degenerating quadrics which was earlier introduced in [H-N-2]. 
   
\noindent{\bf Dedication} 
The question of the behaviour of quadric invariants under degeneration 
was originally posed by M.S. Narasimhan in 1985 in order to
test the rationality of some unirational moduli spaces.
This was motivated by his cohomological criterion for rationality (which he
discovered in the late 1960's), which says that the third integral 
cohomology of a nonsingular complex 
projective variety is torsion free if the variety is rational. 
The attempt to test of this criterion for the Narasimhan-Ramanan 
desingularization of moduli [Na-Ra], 
which has a natural degenerating family of 
conics on it, led to the results described in [N-1] and [N-2]. 
We are happy to dedicate this paper to Professor Narasimhan on his 80th 
birthday, for providing continued encouragement and inspiration over the years.

\section{Purity and Gysin boundary maps}

\centerline{\bf Regular pairs of schemes}

A {\bf regular scheme} 
will as usual mean a locally noetherian scheme all
whose local rings are regular. 
We will work over a regular base scheme $S$ on which $2$ is invertible.   
All schemes, algebraic spaces and  
algebraic stacks and all morphisms will be assumed to be 
over $S$, and the structure morphisms to $S$ will be assumed to be 
quasi-compact. Important extreme examples are when  $S = \spec \Z[1/2]$
or when $S = \spec k$ for an algebraically closed field of 
characteristic $\ne 2$.

A {\bf regular pair $(X,Y)$ of schemes of pure 
codimension $1$} (or simply a `regular pair',
as the codimension will always be assumed to be $1$ in what follows) 
consists of a regular scheme $X$ and a closed subscheme 
$Y$ which  
is regular and of codimension $1$ in $X$ at all points of $Y$.
In particular, $Y$ is an effective Cartier divisor in $X$.
A special class of regular pairs are the {\bf smooth pairs over the base $S$},
where $S$ is a given regular  scheme, and 
$X$ and $Y$ are smooth $S$-schemes with $Y$ a closed subscheme of $X$,
which we will assume to have codimension $1$. An example of 
a regular pair is where $X = \spec \Z[1/2]$ and $Y = \spec \FF_p$ where
$p$ is an odd prime. Note that in this example $X$ is smooth over the base
$S = \spec \Z[1/2]$ but $Y$ is not, so this regular pair is not a smooth 
pair over $S$. Regular or smooth pairs of algebraic spaces have a similar
definition.

Recall that if $(X,Y)$ is a regular pair of algebraic spaces 
of codimension $1$ 
such that $2$ is invertible on $X$, then the absolute 
purity theorem of Gabber (originally conjectured by Grothendieck -- 
see [Fu] for Gabber's proof) gives us the following. 

\stm\label{Absolute purity for schemes}
{\bf Gabber's absolute purity theorem.} 
Under the above hypothesis, 
the local cohomology sheaves $\un{H}^r_Y(\FF_{2,X})$ on $X$ are given by
$$ %
\un{H}^r_Y(\FF_{2,X}) = \left\{ \begin{array}{cc}
0 & \mbox{ for } r \ne 2, \\
i_*\FF_{2,Y} & \mbox{ for } r = 2.
\end{array}
\right.$$
By the resulting degenerateness of the
local-global spectral sequence, this gives the natural identifications
$$H^r_Y(X,\FF_2) = H^{r-2}(X, \un{H}^2_Y(\FF_{2,X})) = H^{r-2}(Y, \FF_2).$$ 
If moreover $Y$ is irreducible, then under the above isomorphism 
the {\bf fundamental class} $s_{Y/X} \in H^2_Y(X,\FF_2)$ of $Y$ in $X$ 
corresponds to the generating section in $H^0(X, \un{H}^2_Y(\FF_{2,X}))$.

\medskip

Under the isomorphism $H^r_Y(X,\FF_2) \cong H^{r-2}(Y, \FF_2)$,
the long exact sequence for cohomology with supports $Y$ becomes
the long exact Gysin sequence 
$$\ldots\to H^r(X,\FF_2) \to H^r(X-Y,\FF_2) \stackrel{\pa}{\to}
H^{r-1}(Y,\FF_2) \to H^{r+1}(X,\FF_2) \to \ldots$$
for the pair $(X,Y)$. 
The maps $\pa : H^r(X-Y,\FF_2) \to H^{r-1}(Y,\FF_2)$
are the {\bf Gysin boundary maps}.

\bigskip

\centerline{\bf Cohomology and local cohomology for algebraic stacks}

We largely follow the terminology of Laumon and Moret-Bailly [L-MB]
for algebraic stacks. Recall that we work over a regular
base scheme $S$ on which $2$ is invertible, and all    
schemes, algebraic spaces and algebraic stacks are assumed to be 
quasi-compact over $S$. The word `stack' will signify such an algebraic stack,
and morphisms will always be over $S$.

Let  $Lis{\mhyphen}Et(\XX)$ denote the smooth-\'etale site of $\XX$, 
and let $\XX_{lis{\mhyphen}et}$ be the resulting topos, which we
will denote simply by $\XX$. 
All sheaves and cohomology for $\XX$ 
will be in the smooth-\'etale topology. 
The sheaf of main interest to us is the constant sheaf 
$\FF_{2,\XX}$ in $\XX_{lis{\mhyphen}et}$.
This sheaf is Cartesian and representable. We will simply denote it by $\FF_2$.
We denote by $Mod_{cart}(\XX,\FF_2)$ the abelian category of Cartesian sheaves
of $\FF_2$-modules on $\XX$, and by $D^+_{cart}(\XX,\FF_2)$ the  
derived category of cohomologically Cartesian bounded below complexes
of sheaves of $\FF_2$-modules on $\XX$.
For an algebraic space $X$, the inclusion functor
$Et(X)\hra Lis{\mhyphen}Et(X)$ from the \'etale to the smooth
site induces a geometric morphism of topoi $X_{lis{\mhyphen}et}\to  X_{et}$
under which we get an equivalence between $X_{et}$ and the subcategory
of Cartesian sheaves in $X_{lis{\mhyphen}et}$ (see [L-MB] 12.2.3).
This induces an equivalence between $D^+(X_{et},\FF_2)$ and 
$D^+_{cart}(X_{lis{\mhyphen}et},\FF_2)$. Consequently, \'etale
cohomology (or hypercohomology) on $X_{et}$ where $X$ is an algebraic space
may be regarded as the smooth-\'etale 
cohomology (or hypercohomology) of the corresponding Cartesian object on 
$X_{lis{\mhyphen}et}$. 

It is known due to Behrend and Gabber that a
$1$-morphism $f: \XX' \to \XX$ of stacks does not functorially 
induce a geometric morphism of topoi from $\XX'_{lis{\mhyphen}et}$ to 
$\XX_{lis{\mhyphen}et}$. However, Olsson showed in [O] how to functorially 
associate to $f$ a pull-back functor
$$f^* :  D^+_{cart}(\XX,\FF_2)\to D^+_{cart}(\XX',\FF_2)$$
of triangulated categories.  
The facts about $f^*$ that we need are only a small part of the 
theory of six operations for algebraic stacks developed by 
Behrend [Be], Olsson [O], Laszlo and Olsson [L-O], and other authors.

{\it Notation:} The functor 
$f^* :  D^+_{cart}(\XX,\FF_2)\to D^+_{cart}(\XX',\FF_2)$ is actually
called $f^{-1}$ in [O]. Our notation $f^*$ for it follows the section
4.3 of [L-O].  

We now briefly recall the construction of $f^*$ from [O]. 
Let $X\to \XX$ be a smooth atlas where $X$ is an
algebraic space over $S$, and let $X^+_{\bullet} = cosk^+_0(X/\XX)$ be its 
nerve regarded as a strictly simplicial algebraic space, 
where $X_n$ is the fiber product of $n+1$ copies of $X$
over $\XX$. 
The word `strict' means that only the 
morphisms $X_n\to X_m$ in $cosk_0(X/\XX)$ which correspond to 
strictly monotonic maps $[m]\to [n]$ (and which are therefore
smooth morphisms) are retained as part of the structure of $X^+_{\bullet}$.
As these morphisms are smooth, the topos
$X^+_{\bullet,lis{\mhyphen}et}$ makes sense. 
Let $X^+_{\bullet, et}$ be the corresponding topos with
\'etale topology.

There are natural geometric morphisms of topoi 
$\pi: X^+_{\bullet,lis{\mhyphen}et} \to \XX_{lis{\mhyphen}et}$ 
and $\epsilon: X^+_{\bullet,lis{\mhyphen}et} \to X^+_{\bullet, et}$.
Using the theory of cohomological descent developed in [SGA 4] and a 
result of Gabber, Olsson proves (Theorem 4.7 of [O]) that 
the induced functors $\pi^*: D^+_{cart}(\XX_{lis{\mhyphen}et}) \to 
D^+_{cart}(X^+_{\bullet,lis{\mhyphen}et})$ and 
$\epsilon^* : D^+_{cart}(X^+_{\bullet, et})\to 
D^+_{cart}(X^+_{\bullet,lis{\mhyphen}et})$
are equivalences of triangulated categories, with quasi-inverses
respectively $R\pi_*$ and $\epsilon_*$. 

{\it Note}: The Cartesian conditions for $\XX_{lis{\mhyphen}et}$ and for   
$X^+_{\bullet, et}$ are meant respectively in the sheaf sense and in 
the strictly simplicial sense, while the Cartesian condition 
for $X^+_{\bullet,lis{\mhyphen}et}$ is meant to be simultaneously satisfied 
in both the sheaf and the strictly simplicial senses.

The composite $\epsilon_*\circ \pi^*$ defines an 
equivalence of triangulated categories
$$\psi: D^+_{cart}(\XX,\FF_2) \to D^+_{cart}(X^+_{\bullet,et},\FF_2).$$
Given a $1$-morphism of $S$-stacks $f: \XX' \to \XX$,  
let $\XX'' = \XX'\times_{\XX}X$ (which is a stack), and let $X' \to \XX''$
be a smooth atlas where $X'$ is an algebraic space. Then we can define a
similar equivalence of triangulated categories
$\psi': D^+_{cart}(\XX',\FF_2) \to D^+_{cart}(X'^+_{\bullet},\FF_2)$.
The induced morphism $f_{\bullet}:  X'^+_{\bullet} \to X^+_{\bullet}$
actually gives a geometric morphism of topoi
$$(f_{\bullet}^*, f_{\bullet,*}) : X'^+_{\bullet, et} \to X^+_{\bullet, et}$$
so it induces a functor of triangulated categories
$$f_{\bullet}^* : D^+_{cart}(X^+_{\bullet,et},\FF_2) \to 
D^+_{cart}(X'^+_{\bullet},\FF_2).$$
Finally, Olsson defines (see 9.16 of [O]) the functor  
$f^* : D^+_{cart}(\XX,\FF_2)\to D^+_{cart}(\XX',\FF_2)$
as the composite
$$D^+_{cart}(\XX,\FF_2)\stackrel{\psi}{\to} D^+_{cart}(X^+_{\bullet, et},\FF_2)
\stackrel{f_{\bullet}^*}{\to} D^+_{cart}(X'^+_{\bullet, et},\FF_2)
\stackrel{\psi'^{-1}}{\to} D^+_{cart}(\XX',\FF_2).$$

It is shown in [O] that the above functor $f^*$ is left adjoint to the functor
$Rf_*: D^+_{cart}(\XX',\FF_2)\to D^+_{cart}(\XX,\FF_2)$. Hence it is 
well-defined, independent of the 
choices of the smooth atlases $X\to \XX$ and $X'\to \XX''$.
Moreover, it is functorial in $f$, that is, given 
$\XX\stackrel{f}{\to} \YY\stackrel{g}{\to}{\mathfrak Z}$, we have 
$(g\circ f)^* = f^*\circ g^*$.

{\noindent{\bf Note}} In the special case where the $1$-morphism 
$f: \XX'\to \XX$ is smooth (but not necessarily representable),
the pullback $f^* : \XX_{lis{\mhyphen}et} \to \XX_{lis{\mhyphen}et}'$ 
on sheaves of sets is indeed (left) exact, and so 
$f^* : D^+_{cart}(\XX,\FF_2)\to D^+_{cart}(\XX',\FF_2)$ 
can be defined directly.
As one would desire, the direct definition of 
$f^*: D^+_{cart}(\XX,\FF_2)\to D^+_{cart}(\XX',\FF_2)$
when $f$ is smooth coincides with the above indirect general definition.

For any $f: \XX'\to \XX$, with the above definition 
of the functor 
$f^*: D^+_{cart}(\XX,\FF_2)\to D^+_{cart}(\XX',\FF_2)$, 
we have $f^*\FF_{2,\XX} = \FF_{2,\XX'}$. 
Hence for any $F$ in $D^+_{cart}(\XX,\FF_2)$, the functor
$f^* : D^+_{cart}(\XX,\FF_2)\to D^+_{cart}(\XX',\FF_2)$ induces a 
pull back homomorphism $f^*$ on hypercohomology 
$${\mathbb H}^i(\XX,F) =   Hom_{D^+_{cart}(\XX,\FF_2)}(\FF_2, F[i]) 
\to   Hom_{D^+_{cart}(\XX',\FF_2)}(\FF_2, f^*F[i])  
= {\mathbb H}^i(\XX',f^*F).$$
In the important special case where $F = \FF_2$, we thus get homomorphisms 
$$f^* : H^i(\XX,\FF_2) \to H^i(\XX',\FF_2).$$

If $i: \YY \hra \XX$ is a closed substack of an algebraic stack $\XX$ over $S$,
and if $F$ is a sheaf in $\XX_{lis{\mhyphen}et}$, then recall that  
the $(0)$-th local cohomology sheaf
$\un{H}^0_{\YY}(F) = i_*i^!(F)$ can be directly defined as follows.
For any $U$ in $Lis{\mhyphen}Et(\XX)$, 
let $V\subset U$ be the pullback of $\YY$.
Then $\un{H}^0_{\YY}(F)(U)$ is defined to be the kernel of the restriction map
$F(U) \to F(U-V)$. The functor $F\mapsto \un{H}^0_{\YY}(F)$ is left exact,
and so we have its derived functor 
$i_*Ri^! : D^+_{cart}(\XX,\FF_2)\to D^+_{cart}(\XX,\FF_2)$. 
({\it Note}: It is enough for our
purpose to directly define 
$i_*Ri^!$, without first defining $Ri^!$ and then composing it with
$i_*$.) 

If $j: \XX-\YY \hra \XX$ is the complementary open inclusion, then for
any $F$ in $D^+_{cart}(\XX,\FF_2)$ we have a functorial exact triangle 
$(i_*Ri^!\FF_{2,\XX} \to \FF_{2,\XX}\to j_*j^*\FF_{2,\XX} \to)$ 
in $D^+_{cart}(\XX,\FF_2)$, which corresponds under the equivalence 
$\psi: D^+_{cart}(\XX,\FF_2)\to D^+_{cart}(X^+_{\bullet, et},\FF_2)$
to the exact triangle 
$(i_{\bullet,*}Ri_{\bullet}^!\FF_{2,X^+_{\bullet,et}} 
\to \FF_{2,X^+_{\bullet,et}} 
\to j_{\bullet,*}j_{\bullet}^*\FF_{2,X^+_{\bullet,et}} \to)$ in 
$D^+_{cart}(X^+_{\bullet, et},\FF_2)$.
Here, $Y\subset X$ is the pullback of $\YY$ under the atlas $X\to \XX$,
$i_{\bullet}: Y^+_{\bullet} \hra X^+_{\bullet}$ is the corresponding 
closed strictly simplicial subspace and 
$j_{\bullet} : X^+_{\bullet}- Y^+_{\bullet} \hra X^+_{\bullet}$ 
is the complementary open inclusion.

The local cohomology groups are defined in terms 
of $i_*Ri^!$ as the hypercohomologies 
$$H^r_{\YY}(\XX,\FF_2) = {\mathbb H}^r(\XX,i_*Ri^!\FF_2).$$ 
Let $i: \YY \hra \XX$ (resp. $i': \YY' \hra \XX'$) be closed 
substacks, such that the closed substack $f^{-1}(\YY)\subset \XX'$ 
has the same support as $\YY'$. We have a natural map 
$$f^*i_*Ri^!\FF_{2,\XX} \to i'_*Ri'^*\FF_{2,\XX'}$$
in $D^+_{cart}(\XX')$ (see below), 
which at the level of hypercohomologies gives a
pull back homomorphism on local cohomologies  
$$f^* : H^i_{\YY}(\XX,\FF_2) \to H^i_{\YY'}(\XX',\FF_2).$$
The homomorphism $f^*i_*Ri^!\FF_{2,\XX} \to i'_*Ri'^*\FF_{2,\XX'}$
has the following strictly simplicial construction.
Choosing atlases as above, 
let $i_{\bullet}: Y^+_{\bullet}\subset X^+_{\bullet}$ and
$i'_{\bullet} : Y'^+_{\bullet}\subset X'^+_{\bullet}$ be the corresponding 
inverse images.
Note that $f$ gives $f_{\bullet} : X'^+_{\bullet}\to X^+_{\bullet}$ under
which the support of $f_{\bullet} ^{-1}(Y^+_{\bullet})$ is 
the support of $Y'^+_{\bullet}$. 
Then $f_{\bullet}$ induces a commutative diagram 
(morphism of exact triangles) in 
$D^+_{cart}(X'^+_{\bullet, et},\FF_2)$
$$%
\begin{array}{cccccc} 
f_{\bullet}^*i_{\bullet,*}Ri_{\bullet}^!\FF_{2,X^+_{\bullet,et}} & \to & 
f_{\bullet}^*\FF_{2,X^+_{\bullet,et}}&\to
     & f_{\bullet}^*j_{\bullet,*}j_{\bullet}^*\FF_{2,X^+_{\bullet,et}} & \to\\
\da &    &\da&       &\da&    \\
i'_{\bullet,*}Ri'^!_{\bullet}\FF_{2,X'^+_{\bullet, et}} & \to & 
\FF_{2,X'^+_{\bullet, et}}&\to& 
j'_{\bullet,*}j'^*_{\bullet}\FF_{2,X'^+_{\bullet, et}} 
& \to 
\end{array}$$
Hence we get the following commutative diagram in 
$D^+_{cart}(\XX',\FF_2)$
$$%
\begin{array}{cccccc} 
f^*i_*Ri^!\FF_{2,\XX} & \to & f^*\FF_{2,\XX}&\to & 
f^*j_*j^*\FF_{2,\XX} & \to\\
\da &    &\da&       &\da&    \\
i'_*Ri'^!\FF_{2,\XX'} & \to & \FF_{2,\XX'}&\to& j'_*j'^*\FF_{2,\XX'} & \to
\end{array}$$

\rem\label{checking commutes} 
To check that a homomorphism in $D^+_{cart}(\XX,\FF_2)$ 
is zero (or to check that a diagram commutes in $D^+_{cart}(\XX,\FF_2)$), 
it is {\it not enough} to check that its pull back 
is zero in $D^+(X_{et},\FF_2)$ for an atlas 
$X\to \XX$ -- rather, one needs to check that 
it is zero in the derived category of
the strictly simplicial topos $X^+_{\bullet, et}$, 
so that cohomological descent can be employed.  
However, if it is a homomorphism (or a diagram) coming from 
$Mod_{cart}(\XX,\FF_2)$, to check that it is zero 
(or commutes) in $D^+_{cart}(\XX,\FF_2)$, 
it is enough to check that its pull back to $Mod(X_{et},\FF_2)$ is zero
(or commutes).

\rem\label{local cohomology commutative diagram} 
By taking hypercohomologies, the above morphism of 
exact triangles gives a commutative diagram
$$%
\begin{array}{ccccccc} 
H^i(\XX,\FF_2)             &\to & H^i(\XX-\YY,\FF_2)                & 
\stackrel{d}{\to} & H^{i+1}_{\YY}(\XX,\FF_2)  &\to& H^{i+1}(\XX,\FF_2)\\
\da &    &\da&       &\da&   & \da\\
H^i(\XX',\FF_2)             &\to & H^i(\XX'-\YY',\FF_2)                & 
\stackrel{d}{\to} & H^{i+1}_{\YY'}(\XX',\FF_2)  &\to& H^{i+1}(\XX',\FF_2)
\end{array}$$

\bigskip

\centerline{\bf Regular pairs of algebraic stacks}

Recall that an algebraic stack $\XX$ is said to be a 
{\bf regular stack} if 
for some (hence for every) smooth atlas $X \to \XX$ where $X$ is
a scheme (or an algebraic space), the total space 
$X$ is regular. A {\bf regular pair
$(\XX,\YY)$ of stacks} will mean a closed embedding $i: \YY \hra \XX$ 
such that both $\XX$ and $\YY$ are regular stacks, and $\YY$ has 
pure codimension $1$. Certain crucial examples of such pairs
for us, namely the pairs $(\M,\D)$ and $(\N,\Zz)$ 
that we introduce later, are in fact
{\bf smooth pairs of $S$-stacks}, that is, the base scheme $S$ 
is regular and the structure morphisms $\XX\to S$ and $\YY\to S$ are smooth.

\stm{\bf (Absolute purity theorem for algebraic stacks.)} 
With all the above preparation, it follows  
from the absolute purity theorem for schemes that  
the absolute purity theorem 
(Statement \ref{Absolute purity for schemes}) holds also for 
regular pairs of stacks $(\XX,\YY)$ of pure codimension $1$
(see [L-O] section 4.10). In particular, 
we get a long exact Gysin sequence as in Statement 
\ref{Absolute purity for schemes}.

\begin{lemma}\label{multiplicity lemma}
Let $(\XX,\YY)$ and $(\XX',\YY')$ be regular pairs of algebraic stacks 
of pure codimension $1$ such that $\YY$ and $\YY'$ are connected, and 
$2$ is invertible on $\XX$ and $\XX'$.
Let $f : \XX' \to \XX$ be a 
$1$-morphism such that the closed substack
$\YY '' = \XX'\times_{\XX} \YY$ of $\YY$ is equal to $m\YY'$ 
for some integer $m \ge 1$
(that is, its ideal sheaf $I_{\YY''}$ equals $I_{\YY'}^m$). 
Then the following statements hold. 

{\noindent (1)}
Under $f^* : H^2_{\YY}(\XX,\FF_2)\to H^2_{\YY'}(\XX',\FF_2)$,
the fundamental classes $s_{\YY/\XX}$ and 
$s_{\YY'/\XX'}$ are related by $s_{\YY/\XX} \mapsto m \cdot s_{\YY'/\XX'}$.

{\noindent (2)} The following diagram commutes
$$%
\begin{array}{ccccccc} 
H^i(\XX,\FF_2)             &\to & H^i(\XX-\YY,\FF_2)                & 
\stackrel{\pa}{\to} & H^{i-1}(\YY,\FF_2)           &\to& H^{i+1}(\XX,\FF_2)\\
~~\da{\scriptstyle f^*}&    & ~~\da{\scriptstyle f_{\XX'-\YY'}^*}&       
                    &~~\da{\scriptstyle mf_{\YY'}^*}&   &
~~\da{\scriptstyle f^*}\\
H^i(\XX',\FF_2)              &\to &   H^i(\XX'-\YY',\FF_2)                & 
\stackrel{\pa'}{\to} & H^{i-1}(\YY',\FF_2)            &\to& H^{i+1}(\XX',\FF_2)
\end{array}$$
where the rows are the respective long exact Gysin sequences. 
In particular, 
$$\pa \circ f_{\XX'-\YY'}^* = m\cdot  f_{\YY}^* \circ \pa'.$$
\end{lemma}

\proof The statement (1) for smooth pairs of schemes
is due to Grothendieck (see [SGA $4{1\over 2}$] page 138), 
who used the purity theorem for smooth pairs. 
If we instead use the absolute purity theorem, the same proof applies
to prove the statement (1) for regular pairs $(X,Y)$ 
of schemes or algebraic spaces.

Now consider the special case when $(X,Y)$ and $(X',Y')$ are
regular pairs of algebraic spaces, and $f: X'\to X$ is a morphism
for which $f^{-1}(Y) = mY'$. As $s_{Y/X} \mapsto m \cdot s_{Y'/X'}$, 
the following diagram commutes, where
the horizontal isomorphisms come from purity, and where 
the right hand vertical map is multiplication by $m$.
$$%
\begin{array}{ccc}
f^*i_*Ri^!\FF_{2,X} & \stackrel{\sim}{\to} & f^*i_*\FF_{2,Y}[-2] \\
\da                 &     & ~ \da {\scriptstyle m}                \\
i'_*Ri'^!\FF_{2,X'} & \stackrel{\sim}{\to} & i_*\FF_{2,Y'}[-2] 
\end{array}$$
Note that the above is actually (the shift by $-2$ of)  
a commutative diagram coming from $Mod_{cart}(X',\FF_2)$. 
Hence it follows by Remark \ref{checking commutes} that the following 
diagram commutes in $D^+_{cart}(\XX',\FF_2)$.
$$%
\begin{array}{ccc}
f^*i_*Ri^!\FF_{2,\XX} & \stackrel{\sim}{\to} & f^*i_*\FF_{2,\YY}[-2] \\
\da                 &     &   ~ \da {\scriptstyle m}              \\
i'_*Ri'^!\FF_{2,\XX'} & \stackrel{\sim}{\to} & i_*\FF_{2,\YY'}[-2] 
\end{array}$$
The statement (1) follows for stacks by applying  
hypercohomology functors ${\mathbb H}^2$ 
to the above diagram and chasing the image of the fundamental class 
$s_{\YY/\XX}$ in the resulting diagram. 

Applying $(i+1)$-th hypercohomology functors
to the above diagram, it follows that the following diagram commutes.
$$%
\begin{array}{ccc}
H^{i+1}_{\YY}(\XX,\FF_2) & \stackrel{\sim}{\to} &  H^{i-1}(\YY,\FF_2) \\
{\scriptstyle f^*}\da~~~&    &~~~\da{\scriptstyle m\, f_{\YY'}^*}\\
H^{i+1}_{\YY'}(\XX',\FF_2) & \stackrel{\sim}{\to} &  H^{i-1}(\YY',\FF_2)
\end{array}$$
By combining the above commutative diagram with 
the one given by Remark \ref{local cohomology commutative diagram}, 
the statement (2) follows.
\hfill$\square$

\rem The conclusions of Lemma \ref{multiplicity lemma}.(2) 
remains true even if in the hypothesis of the lemma the stacks $Y$ and $Y'$ 
are not assumed to be connected and the multiplicity $m$ is
allowed to vary as a locally constant function over $Y'$.
To see this, note
that if $Y_i$ are the connected components of $Y$ and $Y'_j$
are the connected components of $Y'$, and if $Y'_{j_0}$ maps into $Y_{i_0}$,
then
we can apply the lemma to the pairs $(X - \cup_{i\ne i_0}Y_i,\, Y_{i_0})$ and 
$(X' - \cup_{j\ne j_0}Y'_j,\, Y'_{j_0})$, and conclude by excision.

\rem\label{A1 homotopy invariance} 
Any smooth atlas $X\to \XX$ of an $S$-stack functorially gives 
a smooth atlas $X\times \AA^1\to \XX\times \AA^1$ where $\AA^1 = \spec \Z[t]$.
Using these atlases and the corresponding property for algebraic spaces, 
it follows
from the definition of the pullback map on cohomologies of 
stacks that if $f,g \in \Gamma(\XX,\OO_{\XX})$ are regular functions with 
graphs 
$\gamma_f,\gamma_g  : \XX \to \XX\times \AA^1$, then we have
$\gamma_f^* = \gamma_g^* : H^*(\XX\times \AA^1,\FF_2) \to H^*(\XX,\FF_2)$.

\section{Quadric bundles and their degenerations}

\centerline{\bf Quadratic triples and quadric bundles}

Let $S$ be a scheme over $\Z[1/2]$. 
A {\bf quadratic triple} $T = (E,L,b)$ on an $S$-scheme $X$ 
by definition consists of a bilinear form $b$ on a vector bundle 
$E$ on $X$, taking values in a line bundle $L$ on $X$, 
such that $b$ is nowhere zero on $X$, that is, 
$b_x$ is not identically zero on any fiber $E_x$. 
We call $\rank(E)$ (which we assume to be globally constant) the 
{\bf dimension} of the triple $T$ and the $\rank(b_x)$ of the quadratic
form $b_x$ on the fiber $E_x = k(x)\otimes_{\OO_X}E$ as the {\bf rank} 
of $T$ at $x \in X$. 
We will say that a quadratic triple $T = (E,L,b)$ on $X$ is 
{\bf nondegenerate} if at each $x\in X$ we have
$\rank(b_x) = \rank(E)$, equivalently, the $\OO_X$-linear map
$b: E \to L\otimes_{\OO_X}E^*$ is an isomorphism.

Given a triple $T =(E,L,b)$ and a line bundle $K$ on $X$, the triple  
$T\otimes K$ is defined to be 
$(E\otimes K,L\otimes K^{\otimes 2}, b\otimes 1_K )$. 
We say that two 
triples $T = (E,L,b)$ and $T' = (E',L',b')$ on $X$ are 
{\bf equivalent triples} if
there exists a line bundle $K$ on $X$ 
such that $T\otimes K$ is isomorphic to $T'$. 
A {\bf quadric bundle} $Q$ over $X$ is by definition an 
equivalence class $[E,L,b]$ of triples such that 
for each $x\in X$ we have $b_x \ne 0$.
As $b_x$ is not identically zero on any fiber $E_x$ by assumption, 
we see that $[E,L,b]$ defines a closed subscheme
$C \subset P(E) = \Proj \Sym (E^*)$ by 
the vanishing of the quadratic form corresponding to $b$,
which is flat over $X$ and whose each fiber is a 
quadric hypersurface $C_x\subset P(E_x)$. Conversely
the equivalence class $Q=[E,L,b]$ can be recovered from 
a closed subscheme $C \subset P(E)$ in a banal projective bundle,
such that $C$ is flat over $X$ and every schematic 
fiber is a quadric hypersurface.

\noindent{\bf Quadratic triples and quadric bundles over algebraic stacks:} 
In the above definitions, the base space for quadratic triples
and quadric bundles 
was taken to be an $S$-scheme $X$. But the same definitions work equally
well when $X$ is an algebraic stack over $S$. In what follows, the base
spaces for our quadratic triples or quadric bundles will be assumed to
be algebraic stacks over $S$, unless otherwise indicated.

\bigskip

\centerline{\bf Minimally degenerate quadric bundles,  
associated nondegenerate triples}

The {\bf discriminant $\det(T)$} of a quadratic triple $T =(E,L,b)$ on $X$
is defined by 
$$\det(T) = \det(b) \in \Gamma(X, L^{\rank(E)} \otimes \det(E)^{-2}).$$
A quadratic triple $T$ is nondegenerate if and only if $\det(T)$ 
is nowhere vanishing. Note that 
$$\det(T) = \det(T\otimes K)$$
for any line bundle $K$ on $X$, so for any quadric bundle 
$Q = [E,L,b]$. This enables us to define {\bf the discriminant $\det(Q)$
of a quadric bundle} to be $\det(T)$ for any representative $T$.

A quadratic triple $(E,L,b)$ on base $Y$ is called a 
{\bf minimally degenerate triple} if (i) 
$\det(b) = 0 \in \Gamma(Y, L^{\rank(E)} \otimes \det(E)^{-2})$
and (ii) $\rank(b_y) = n-1$ for all $y\in Y$. (Note that if $Y$ is reduced
then (ii) implies (i), but not otherwise.) 
The equivalence class 
$Q = [E,L,b]$ of such a triple 
will be called a {\bf minimally degenerate quadric bundle}
(this is well defined, independent of the choice of a representative
$(E,L,b)$ of $Q$).
If we regard $b$ as an $\OO_Y$-linear homomorphism
$b: E \to L\otimes E^*$, then $\ker(b)$ is a line bundle, 
the quotient $\ov{E} = E/\ker(b)$ is a vector bundle of rank equal to 
$\rank(E)-1$, and $b$ induces an $L$-valued
quadratic form $\ov{b}$ on $\ov{E}$, giving a  
nondegenerate quadratic triple 
$(\ov{E}, L, \ov{b})$.

The next lemmas follows from a simple computation which we omit.

\begin{lemma}
If $Q = [E,L,b]$ is a minimally degenerate quadric bundle
on $Y$, then the nondegenerate quadratic triple 
$${\mathbb T}^Q = (\ov{E}, L, \ov{b})\otimes \ker(b)^{-1}$$
is well defined, that is, it remains 
unaltered if the representative
$(E,L,b)$ is replaced by $(E,L,b)\otimes K$ for any line bundle $K$ on $Y$.
\hfill $\square$
\end{lemma} 

We will call the above quadratic triple  
${\mathbb T}^Q = (\ov{E}, L, \ov{b})\otimes \ker(b)^{-1}$
(which was originally introduced in [H-N-2]) 
the {\bf associated nondegenerate quadratic triple} to 
the minimally degenerate quadric bundle $Q$.

\bigskip

\centerline{\bf Mildly degenerating triples on a regular pair $(X,Y)$}

A {\bf mildly degenerating triple} $(E,L,b)$ on a
regular pair $(X,Y)$ of stacks is a 
quadratic triple $T =(E,L,b)$ on $X$ such that $T|_{X-Y}$ is nondegenerate
while $T|_Y$ is minimally degenerate, that is,  
$\rank_x(b) = \rank(E)$ for each $x\in X-Y$ while 
$\rank_y(b) = \rank(E) - 1$ for each $y\in Y$. As such a 
$Y$ is reduced, the discriminant $\det(T)$ vanishes on $Y$.

Given a mildly degenerating 
triple $T = (E,L,b)$ on a regular pair $(X,Y)$ of schemes, 
the vanishing
multiplicity of the discriminant $\det(T)$ along the divisor 
$Y$ defines a 
function $\nu_Y(T) : |Y| \to \Z_{\ge 0}$ on the 
set $|Y|$ of connected components of $Y$.
This function will be called the {\bf degeneration multiplicity} of the triple.
If $(X,Y)$ is a regular pair and $m: |Y| \to \Z_{\ge 0}$ 
is a continuous function, then we get a new ideal sheaf $I_Y^m\subset \OO_X$,
and we denote by $mY\subset X$ the corresponding closed subscheme. 
Note that if $m=0$ then $mY = \emptyset$ and if $m=1$ then $mY =Y$.
The above has an obvious interpretation also for mildly degenerating 
triple on a regular pair of stacks.

\rem In particular, if $(X,Y)$ is a regular pair and 
if $Q$ is a mildly degenerating quadric bundle on $(X,Y)$, 
then the restriction $Q_Y$ is a minimally degenerate quadric bundle
on $Y$, which gives us an associated nondegenerate quadratic triple 
${\mathbb T}^{Q_Y}$.

\section{Cohomological invariants}

\centerline{\bf The group schemes $GO_n$ and stacks $BGO_n$}

The {\bf general orthogonal group-scheme} 
$GO_n$ over $\Z[1/2]$ is the closed sub-group scheme of $GL_n$ over 
$\Z[1/2]$, whose 
$R$-valued points for any $\Z[1/2]$-algebra $R$ 
(that is, a commutative ring $R$ in which $2$ is invertible)
are all the $n\times n$
matrices $g$ over $R$ such that ${^t}gg = \lambda I$ for some 
invertible element $\lambda \in R^{\times}$. 
We have a surjective homomorphism $\sigma: GO_n\to \GG_m$ of group schemes 
defined on valued points by $g\mapsto \lambda$ where ${^t}gg = \lambda I$.
This has kernel $O_n$, the orthogonal group scheme over $\Z[1/2]$. 
As both $O_n$ and $\GG_m$ are smooth over $\Z[1/2]$, it follows that 
$GO_n$ is also smooth over $\Z[1/2]$. In fact, it can be seen that $GO_n$ is
a reductive group scheme over $\Z[1/2]$. 

For any scheme $S$ over $\Z[1/2]$, the group-scheme $GO_{n,S}$ 
is obtained by base change, and has a similar description in terms of 
$R$-valued points for affine schemes $\spec R\to S$. 
The group schemes $GL_{n,S}$, $\GG_{m,S}$, $GO_{n,S}$, $O_{n,S}$, etc. 
will be simply denoted by $GL_n$, $\GG_m$, $GO_n$, $O_n$, etc.
for simplicity of notation.

When $n = 2m+1$ is odd, the homomorphism
$\GG_m \times SO_{2m+1} \to GO_{2m+1}$ defined in terms of
valued points by $(\lambda, g) \mapsto \lambda g$ is an isomorphism. 
This induces a $1$-isomorphism of algebraic stacks 
$B\GG_m \times BSO_{2m+1} \to BGO_{2m+1}$. 

If $X$ is a scheme over $S$, then a principal 
$GO_n$-bundle $P$ on $X$, locally trivial 
in the \'etale topology, is equivalent to a nondegenerate quadratic triple
$(E,L,b)$. To any such $P$, 
we functorially associate the vector bundle $E$ given by the defining 
representation $GO_n\hra GL_n$, and $L$ given by the character 
$\sigma : GO_n\to \GG_m$ defined above. The standard
bilinear form $\sum x_iy_i$ defines $b: E\otimes E \to L$. 
This defines a $1$-isomorphism of stacks from the algebraic stack $BGO_n$ to
the algebraic stack $\M -\D$ of nondegenerate quadratic triples
that is introduced later. The inverse $1$-morphism
has an obvious construction using the Gram-Schmidt orthogonalization process. 
The principal bundle $P$ associated to
any $(E,L,b)$ is \'etale locally trivial as the process of
Gram-Schmidt orthogonalization requires taking
square-roots of invertible regular functions, 
which is possible \'etale locally as $2$ is invertible over $X$.

Similarly, the stack $BO_n$ can be identified with the stack
of {\bf nondegenerate quadratic pairs}, where such a pair $(E,q)$
on an $S$-scheme $X$ is just a nondegenerate triple
$(E,\OO_X,q)$.   

On the algebraic stack $BGO_n$, we have a universal nondegenerate quadratic 
triple $(E_n,L_n,b_n)$. Let $L_n$ also denote the total space of the 
line bundle $L_n$, so that $L_n$ is again an algebraic stack with a projection
$\pi : L_n \to BGO_n$ and a zero section $BGO_n \to L_n$, which are 
$1$-morphisms of algebraic stacks. The stack $L_n$ can be obtained from $BGO_n$ 
by a local construction in the sense of Laumon, Moret-Bailly [L-MB] Chapter 14.

Let $L_n^{\times} \subset L_n$ be the open substack which is the complement 
of the zero section $BGO_n \hra L_n$. For any $X$ over $S$, an object
of $L_n^{\times}(X)$ is a tuple $(E,L,b,s)$ where $(E,L,b)$ is a 
nondegenerate quadratic triple on $X$, and $s \in \Gamma(X,L)$
is a nowhere vanishing section. Equivalently, $s : \OO_X \to L$
is an isomorphism. Therefore there is a $1$-morphism of stacks
$BO_n \to L_n^{\times}$, which associates to a nondegenerate
quadratic pair $(F,q: F\otimes F \to \OO_X)$ on $X$ the 
tuple $(F,\OO_X, q, \id_{\OO_X})$. In the reverse direction, 
we can associate to a tuple $(E,L,b,s: \OO_X \to L)$ the 
quadratic pair $(E, s^{-1}\circ b)$, showing that the $1$-morphism
$BO_n \to L_n^{\times}$ is an isomorphism of algebraic stacks. Using this,
we identify $L_n^{\times}$ with $BO_n$. 

As $GO_n$ is a smooth group scheme over $S$, the stack $BGO_n$ 
is smooth over $S$, hence the total space $L_n$ is also smooth over 
$S$. We regard $BGO_n\hra L_n$ as the zero section, so that
$(L_n, BGO_n)$ is a smooth pair over
$S$, in particular, it is a regular pair.
The projection $\pi: L_n \to BGO_n$ induces an isomorphism
$ H^*(BGO_n,\FF_2) \to H^*(L_n,\FF_2)$.
Identifying 
$L_n^{\times}$ with $BO_n$, 
the Gysin sequence for $(L_n, BGO_n)$
becomes the following exact sequence.
{\small 
$$ \ldots \to H^*(BGO_n,\FF_2)\to  
H^*(BO_n,\FF_2) \stackrel{d_n}{\to} 
H^{*-1}(BGO_n,\FF_2) \to H^{*+1}(BGO_n,\FF_2) \to \ldots $$
}

\bigskip

\bigskip

\centerline{\bf Characteristic classes for 
nondegenerate triples and $H^*(BGO_n)$}

Given a group-scheme $G$ over a base $S$,
let $|BG|$ denote the 
set-valued functor on $(QuasiProj_S)^{opp}$
which associates to any quasi-projective $S$-scheme $X$ the 
set of all isomorphism classes of principal $G$-bundles on $X$.
Assuming that $2$ is invertible on $S$,
a {\bf functorial characteristic class} 
for $G$ will mean a natural transformation
from $|BG|$ to the functor $H^*(-,\FF_2)$ which associates
to any $X$ its \'etale cohomology $H^*(X,\FF_2)$. All 
such functorial characteristic
classes form a graded ring $CharClass_G^*$ under the ring operations on 
the various $H^*(X,\FF_2)$. 

If $G$ is of finite-type and flat over $S$, then $BG$ is an algebraic stack 
over $S$. Let $H^*(BG,\FF_2)$ be its cohomology in the 
smooth-\'etale topology. We have a natural homomorphism
$$H^*(BG,\FF_2) \to CharClass_G^*$$
of graded rings, 
which associates to any $\alpha \in H^*(BG,\FF_2)$ the functorial 
characteristic class 
$\alpha'$ defined by $\alpha'(P) = \chi_P^*(\alpha)$ where $\chi_P : X\to BG$
is the classifying morphism of a $G$-bundle $P$ on $X$.

When the base $S$ is of the form $\spec k$ for an algebraically
closed field of characteristic $\ne 2$ and when $G$ is reductive over $k$, 
it was shown in [Bh] that the homomorphism 
$H^*(BG,\FF_2) \to CharClass_G^*$ is an isomorphism. 
We now revisit the argument in [Bh], to show it can be made to work 
over any excellent Dedekind domain on which
$2$ is invertible (for example, $\Z[1/2]$). 
This gives us the following.

\begin{proposition}\label{functorial char classes} 
Let $G$ be a reductive group scheme over a base which is a field or 
an excellent Dedekind domain on which $2$ is invertible. 
Then the natural homomorphism $H^*(BG,\FF_2) \to CharClass_G^*$
is an isomorphism of graded rings.
\end{proposition}

\proof (Sketch) The following three comments (1), (2) and (3)
explain how the proof of the 
corresponding result (Proposition 3.1) in [Bh] can be modified so that
it applies in our case.  

\noindent{(1)} The Lemma 3.1 of [Bh] says that any closed reduced subscheme 
$Z\subset \AA^n_k$
has a finite filtration by locally closed reduced subschemes
$Z = Z_0\supset Z_1\supset \ldots \supset Z_{\ell} = \emptyset$ such that each 
$Z_i-Z_{i+1}$
is smooth over $k$. A modified statement holds  
for closed reduced subschemes $Z\subset \AA^n_{S}$, 
with the condition `$Z_i-Z_{i+1}$ is smooth over $k$' replaced  
by the condition that $Z_i-Z_{i+1}$ is a regular scheme 
(we do not require $Z_i-Z_{i+1}$ to be smooth over $S$). 
Such a filtration exists because in any excellent scheme, 
regular points form an open subscheme (see Grothendieck [EGA-IV$_2$] 
page 215 statement 7.8.3.(iv)). An application of the 
absolute purity theorem of Gabber now completes the argument.

\noindent{(2)} The use of GIT for quotients for linear actions
of reductive groups (in which which [Bh] follows Totaro [T])
is replaced by the use of Seshadri's version 
of GIT in [Se], which in particular works over the base $S$.

\noindent{(3)} Jouanolou's trick   
(in which [Bh] again follows Totaro [T]) remains valid 
over the base $S$ (see [Ju]). 

With the above modifications, the rest of the argument now 
works as in [Bh]. \qed

\bigskip

\bigskip

\bigskip

\bigskip

\centerline{\bf Cohomological invariants for nondegenerate quadric bundles}

Tensor product of triples and line bundles defines a 
$1$-morphism of algebraic stacks
$$\mu: BGO_n\times B\GG_m \to BGO_n$$
which sends $(T,K)\mapsto T\otimes K$ where $T$ is 
a quadratic triple on $X$ and $K$ is a line bundle on $X$.
Following Toda [Td] as in [H-N-2], we say that an element 
$\alpha\in H^*(BGO_n, \FF_2)$ is a {\bf primitive class} if 
under the cohomology homomorphisms induced by $\mu$ and by
the projection $p_1: BGO_n\times B\GG_m \to BGO_n$, we have
$$\mu^*(\alpha) = p_1^*(\alpha).$$
The primitive classes form a subring 
$PH^*(BGO_n, \FF_2)\subset H^*(BGO_n, \FF_2)$.
Over a field or an excellent Dedekind domain $S$ on which $2$ is invertible, 
the arguments used for proving that the natural homomorphism 
$H^*(BG,\FF_2) \to CharClass_G^*$ is an isomorphism can be used to see that 
elements $\alpha \in  PH^*(BGO_n, \FF_2)$ can also be 
characterized as natural transformations from 
$|BGO_n|$ to $H^*(-,\FF_2)$ on $QuasiProj_{S}$ 
such that 
$\alpha(T\otimes K) = \alpha(T)$ for any nondegenerate triple $T$ 
and a line bundle $K$ on a quasi-projective scheme $X$ over $S$.

\section{The main theorem}

\subsection{Statement of the main theorem}

We are now able to state the main theorem, 
having made all the above preparation over \'etale
cohomology and Gysin sequences for algebraic stacks. 
But the statement
is just the algebraic stacky analogue of the corresponding topological  
theorem in [H-N-2], where the complex Lie groups $GO_n(\C)$, etc. are replaced
by the corresponding group schemes over $S$, and singular cohomology 
with coefficients $\FF_2$ is replaced by \'etale cohomology 
with coefficients $\FF_2$.
We will adhere to the notation in [H-N-2] as much as possible.

Let $PH^*(BGO_n, \FF_2)\subset H^*(BGO_n, \FF_2)$ be the 
subring of primitive classes. Its elements are the 
universal cohomological invariants
for nondegenerate quadric bundles of rank $n$.
Let $B(v)^* : H^*(BGO_n) \to  H^*(BO_{n-1})$ be the ring
homomorphism induced by the group homomorphism $v: O_{n-1} \to GO_n$
defined by $g\mapsto \left(\begin{array}{cc} 1& \\ & g \end{array}\right)$.
Let $(E_{n-1},L_{n-1},b_{n-1})$ be the universal triple over $BGO_{n-1}$.
As shown above, the 
complement $L_{n-1}^{\times}$ of the zero section of $L_{n-1}$ is isomorphic to
$BO_{n-1}$, and there is a Gysin boundary map
$d_{n-1} : H^*(BO_{n-1}) \to H^{*-1}(BGO_{n-1})$. 
Finally, let $\delta : PH^*(BGO_n) \to H^{*-1}(BGO_{n-1})$ be the
composite linear map 
$$PH^*(BGO_n) \hra H^*(BGO_n) \stackrel{B(v)^*}{\to} H^*(BO_{n-1})
\stackrel{d_{n-1}}{\to} H^{*-1}(BGO_{n-1})$$ 
With the above notations, we have the following,
where all cohomologies are smooth \'etale cohomologies 
with coefficients $\FF_2$.

\begin{theorem}\label{main theorem} ({\bf Main Theorem}) 
Let $S$ be a regular scheme on which $2$ is invertible.
Let $(X,Y)$ be a regular pair of stacks over $S$. 
Let $Q$ be a quadric bundle on $X$, nondegenerate on $X-Y$
of rank $n \ge 2$, which is minimally degenerate on $Y\subset X$.
Let ${\mathbb T}^{Q_Y}$ be the associated rank $n-1$ nondegenerate 
quadratic triple on $Y$. 
Let $\alpha\in PH^*(BGO_n)$ be a universal quadric invariant, and let 
$\alpha(Q_{X-Y})\in H^*(X-Y)$ be its value on $Q_{X-Y}$.
Let $\delta : PH^*(BGO_n)\to H^{*-1}(BGO_{n-1})$ 
be the linear map defined above,
and let $(\delta(\alpha))({\mathbb T}^{Q_Y})$ be the value of the resulting
$GO_{n-1}$-characteristic class $\delta(\alpha)$ on 
${\mathbb T}^{Q_Y}$. Let $\nu_Y(\det(Q)) \in H^0(Y)$ be the  
vanishing multiplicity along $Y$ of the discriminant 
$\det(Q) \in \Gamma(X,L^n\otimes \det(E)^{-2})$. 
Then under the Gysin boundary map $\pa : H^*(X-Y) \to H^{*-1}(Y)$, 
we have the equality
$$\pa(\alpha(Q_{X-Y})) = 
\nu_Y(\det(Q)) \cdot (\delta(\alpha))({\mathbb T}^{Q_Y})$$
\end{theorem}

The rest of this paper is devoted to proving the theorem.

\subsection{The algebraic stack of quadratic triples}

Let $n\ge 1$. For any $S$-scheme $X$, let $\M_n(X)$ be the groupoid
whose objects are all quadratic triples $T=(E,L,b)$ on $X$ with
$\rank(E) =n$, and a morphism  $T \to T'$, where 
$T = (E,L,b)$ and $T' = (E',L',b')$, is pair of
$\OO_X$-linear isomorphisms $E\to E'$ and $L \to L'$ 
which takes $b$ to $b'$. 
If $f: X' \to X$ is a morphism of $S$-schemes 
and $T$ is a quadratic triple over $X$, the pull-back triple 
$f^*(T) = (f^*E,f^*L, f^*b)$ is defined in the obvious manner. 
For each $f: X'\to X$ and $f'' : X'' \to X'$ over $S$, we can define an  
isomorphism $\phi_{f,f'} : (ff')^*(T) \to f'^*f^*(T)$, which 
makes $\M_n$ into a groupoid-valued pseudo functor on $S$-schemes.
It is clear that the pull-back $f^*$ preserves
both the dimension (rank of $E$) and the point-wise ranks
(ranks of $b_x$), and pulls back the discriminant, so properties 
such as nondegenerateness, minimal degenerateness, etc. are preserved 
under pullbacks. 
The {\bf almost nondegenerate triples}, which
will mean triples whose point-wise rank is $\ge n-1$,  
form a full subgroupoid $\M$ of $\M_n$.
The $S$-groupoid $\M$ has a full subgroupoid $\D$ formed by all the 
quadratic triples $T = (E,L,b)$ on $S$-schemes $X$ such that
$\det(T) =0 \in \Gamma(X, L^n\otimes \det(E)^{-2})$ and 
$\rank(b_x) = n-1$ for all $x\in X$. This is the groupoid of  
minimally degenerate triples.

\begin{proposition}\label{stacks M and D via local construction}
For any $n\ge 1$, almost nondegenerate triples of dimension $n$
form an algebraic stack $\M$ over $S$, and minimally degenerate
triples of dimension $n$ form a closed substack $\D\subset \M$
defined by the vanishing of the discriminant of the universal triple
$\T$ on $\M$. Moreover, $\D$ is a smooth relative divisor in $\M$ over $S$.
The stack $\M$ is a global quotient stack under the reductive group scheme
$GL_n\times \GG_m$ over $S$.
\end{proposition}

\noindent{\bf Proof } We will see below 
that the $S$-groupoid $\M_n$ of all triples is an algebraic stack on 
$S$. It has a filtration by closed algebraic substacks 
$\M_n = \M_n^{\le n} \supset \M_n^{\le n-1} \supset \M_n^{\le n-2} \supset \ldots$
where for any $X$, $\M_n^{\le r}(X)$ is the groupoid of all triples 
$(E,L,b)$ on $X$ such that the determinants of minors of $b$ of size $r+1$ 
vanish
identically on $X$ (the relevant Fitting ideal is $0$). 
The almost nondegenerate
triples form the open substack 
$$\M = \M_n - \M_n^{\le n-2}.$$
The closed substack $\D \subset \M$ is given by 
$$\D = \M_n^{\le n-1} - \M_n^{\le n-2}.$$ 
We now construct these stacks $\M_n^{\le r}$. 
Let $E$ be the universal vector bundle
on the stack $BGL_n$ and let $L$ be the universal line bundle
on $B\GG_m$. On the stack $BGL_n\times B\GG_m$ we get
the vector bundles $p_1^*\un{Sym}^2(E)$ and $p_2^*L$. Then $\M_n$ is
the `total space' (rather,`total stack') 
of the vector bundle corresponding to the locally free 
$\OO$-module $\unhom(p_1^*\un{Sym}^2(E),\,p_2^*L)$ on 
$BGL_n\times B\GG_m$, which is an algebraic stack as it is a 
local construction in the sense of Laumon, Moret-Bailly [L-MB] Chapter 14.
The closed substack $\M_n^{\le r}(X)$ is defined by the 
vanishing of all determinants of minors of size $r+1$. In particular, 
$\D$ is defined in $\M$ by the vanishing of the discriminant, so
it is a principal divisor. 

As $BGL_n\times B\GG_m$ is smooth over $S$, and as
$\M_n$ is a geometric vector bundle over it, $\M_n$ is a smooth stack over $S$.
The Jacobian criterion applied to the discriminant shows that 
$\D$ is a smooth relative divisor in $\M$ over $S$. 

The affine space $\AA^{n(n+1)/2}_S$ over $S$ of 
$n\times n$-symmetric matrices has an action of $GL_n\times \GG_m$
given in terms of valued points by $b \cdot (g,\lambda) = {^t}gbg\lambda^{-1}$.
If $M\subset \AA^{n(n+1)/2}_S$ is the open subscheme where $\rank(b)\ge n-1$,
and if $D\subset M$ is the divisor
defined by $\det(b) =0$, then it can be seen that $\M$ is isomorphic to 
the quotient stack $[M/GL_n\times \GG_m]$ under which $\D$
becomes its the closed substack $[D/GL_n\times \GG_m]$.
\qed

Note that we have a tautological quadratic triple 
$\T_n = (\E_n,\L_n,\beta_n)$ on $(\M,\D)$. 
If $X$ is an algebraic stack, then to any almost nondegenerate 
quadratic triple $T$ on $X$ there functorially corresponds its 
{\bf classifying morphism} $\chi_T : X\to \M$, together
with an isomorphism of triples $\chi_T^*(\T_n) \to T$. 
The degeneration multiplicity along $Y$ of a mildly degenerating 
triple $T$ on a regular pair $(X,Y)$ of stacks has the following obvious
interpretation in terms of its classifying morphism $\chi_{T} : X\to \M$.

\begin{proposition} Let $T = (E,L,b)$ be a mildly degenerating 
triple on a regular pair $(X,Y)$, with degeneration multiplicity
$\nu_Y(T): |Y| \to \Z_{\ge 0}$. Then under the classifying morphism 
$\chi_T : X\to \M$, the pull back of 
the closed substack $\D\subset \M$
is the closed substack $\nu_Y(T)Y\subset X$ defined by
the ideal sheaf $(I_Y)^{\nu_Y(T)}\subset \OO_X$.  \hfill $\square$
\end{proposition}

We now translate the main theorem in terms of the Gysin homomorphism 
$\pa : H^*(\M-\D) \to H^{*-1}(\D)$. For this, we first define 
a $1$-morphism $\rho : \D \to BGO_{n-1}$ as follows. 
If $Y$ is any (affine) scheme
over $S$, an object of the groupoid $\D(Y)$ is a minimally degenerate triple
$T = (F,L,q)$ over $Y$, with $\rank (F) =n$ and $rank(q) = n-1$.
Note that $\ker(q)\subset F$ is a line subbundle, and $q$ induces a 
nondegenerate $L$-valued form $\ov{q}$ on the quotient $F/\ker(q)$.
We get a new nondegenerate triple of rank $n-1$
$$\rho(F,L,q) = (F/\ker(q), L, \ov{q})\otimes \ker(q)^{-1}$$
over $Y$. This is functorial in the triples and respects their pull-backs,
so it defines $1$-morphism of stacks $\rho : \D \to BGO_{n-1}$.
Note that if $Q= [T]$ is a mildly degenerating quadric bundle on $(X,Y)$, 
then  
$${\mathbb T}^{Q_Y} = \rho(T|_Y).$$
Hence the $1$-morphism $\rho : \D \to BGO_{n-1}$ is the 
classifying morphism of the triple ${\mathbb T}^{\Qq_{\D}}$ on $\D$, 
where $\Qq$ is the  mildly degenerating quadric bundle 
$\Qq = [\T]$ on $(\M,\D)$ corresponding to
the universal triple $\T$ on $\M$.

We will prove the main theorem
in the following equivalent re-formulation.

\begin{theorem}\label{universal gysin form of main theorem}
{\bf Main theorem restated in terms of the pair $(\M,\D)$ :}\\
The Gysin boundary map $\pa : H^*(\M -\D) \to H^{*-1}(\D)$ 
on the primitive cohomology $PH^*(\M -\D) \subset H^*(\M -\D)$
is given by the formula
$$\pa|_{PH^*(\M -\D)} = \rho^* \circ d_{n-1} \circ B(v)^*$$
where $B(v)^* : H^*(\M -\D) = H^*(BGO_n) \to H^*(BO_{n-1})$
is induced by $v: O_{n-1} \to GO_n : g \mapsto 
\left(\begin{array}{cc} 1& \\ & g \end{array}\right)$, 
the map $d_{n-1} : H^*(BO_{n-1})\to H^{*-1}(BGO_{n-1})$ is the 
Gysin boundary map for the pair $(BO_{n-1}, BGO_{n-1})$, and 
$\rho :  \D \to BGO_{n-1}$ is the $1$-morphism defined by 
$\rho(F,L,q) = (F/\ker(q), L, \ov{q})\otimes \ker(q)^{-1}$.
\end{theorem}

\noindent{\bf Proof of equivalence:} Consider the universal quadratic triple 
$\T$ on $\M$, and the resulting quadric bundle $\Qq = [\T]$, which 
is a mildly degenerating quadric bundle on the regular pair of stacks 
$(\M,\D)$, with degeneration multiplicity $\nu_{\D}(\det(\Qq))=1$. 
Let $\alpha \in PH^*(\M -\D) = PH^*(BGO_n)$. By definition, 
$\alpha(\Qq_{\M -\D}) = \alpha$. 
Therefore the formula 
$\pa(\alpha) = \nu_Y(\det(Q)) \cdot (\delta(\alpha))({\mathbb T}^{Q_Y})$
of the Theorem \ref{main theorem} applied to $\alpha$ and $\Qq$ becomes
$\pa(\alpha(\Qq_{\M -\D})) = \delta(\alpha)({\mathbb T}^{\Qq_{\D}})$.
As remarked above, ${\mathbb T}^{\Qq_{\D}} = \rho(\T|_{\D})$.
Hence we get for all $\alpha \in PH^*(\M -\D)$ the equality
$\pa(\alpha) = \rho^*\delta(\alpha) = \rho^* \circ d_{n-1} \circ B(v)^*(\alpha)$
by definition of the map $\delta$ from Theorem \ref{main theorem}.
As this holds for all $\alpha \in PH^*(\M -\D)$, we get
$\pa|_{PH^*(\M -\D)} = \rho^* \circ d_{n-1} \circ B(v)^*$ as desired.
(The above is a typical `Yoneda argument'.)

Conversely, suppose the above formula holds. Given any mildly
degenerating quadric $Q$ on a regular pair of stacks $(X,Y)$, 
let $\chi : X \to \M$ be the classifying $1$-morphism for any chosen
quadratic triple $T$ on $X$ which represents $Q$.
Let $m = \nu_Y(Q)$. 
By Lemma \ref{multiplicity lemma}.(2) we get
a commutative diagram
$$\begin{array}{ccc}
H^i(\M-\D)&\stackrel{\pa_{(\M,\D)}}{\to} & H^{i-1}(\D)\\
{\scriptstyle \chi|_{X-Y}^*}\da ~~~&
                             &~~~\da{\scriptstyle m\cdot\chi|_Y^*}\\
H^i(X-Y)&\stackrel{\pa_{(X,Y)}}{\to}&H^{i-1}(Y)
\end{array}
$$
We now make the following observation (\ref{pulling back cohomology})
before continuing with the proof that (\ref{main theorem}) 
implies (\ref{universal gysin form of main theorem}).
\stm\label{pulling back cohomology} 
Note that for any $\alpha \in H^i(\M-\D)$, by definition 
$\alpha(T) = \chi|_{X -Y}^*(\alpha)$, and in particular when 
$\alpha \in PH^i(\M-\D)$, by definition 
$\alpha(Q) = \chi|_{X-Y}^*(\alpha)$. Also, for any 
$\beta \in H^{i-1}(\D)$, 
by definition $\beta(T|_Y) = \chi|_Y^* (\beta)$.

\medskip

For any $\alpha \in PH^*(\M-\D)$ 
we therefore have the following sequence of equalities.
\begin{eqnarray*}
\pa_{(X,Y)}(\alpha(Q)) 
& = & \pa_{(X,Y)}\chi|_{X -Y}^*(\alpha) 
\mbox{ by the Statement \ref{pulling back cohomology},} \\
& = & m\cdot \chi|_Y^*\pa_{(\M,\D)}(\alpha) 
 \mbox{ by the above application of Lemma \ref{multiplicity lemma}.(2),}\\
& = & m\cdot \chi|_Y^*\pa|_{PH^*(\M-\D)}(\alpha) 
     \mbox{ as by assumption } \alpha \in PH^*(\M-\D), \\
& = & m\cdot \chi|_Y^*\rho^* \circ d_{n-1} \circ B(v)^*(\alpha)
      \mbox{ if Theorem \ref{universal gysin form of main theorem} holds,}\\
& = & m\cdot \chi|_Y^*\rho^* \delta (\alpha),\\
&   & \mbox{ where }\delta : PH^*(\M-\D) = PH^*(BGO_n) \to H^{*-1}(BGO_{n-1})\\
&   & \mbox{ is defined as earlier by } \delta = d_{n-1} \circ B(v)^* .
\end{eqnarray*}
This gives the formula 
$$\pa_{(X,Y)}(\alpha(Q)) = m\cdot \chi|_Y^*\rho^* \delta (\alpha)$$
Note that the composite $1$-morphism 
$Y \stackrel{\chi|_Y}{\to} \D  \stackrel{\rho}{\to} BGO_{n-1}$
is exactly the characteristic morphism $\chi' : Y\to BGO_{n-1}$
of the reduced triple 
${\mathbb T}^{Q_Y} = (E/\ker(b), L, \ov{b})\otimes \ker(b)^{-1}$
where $T = (E,L,b)$ represents $Q$. 
Hence we get 
$$\chi|_Y^*\rho^* \delta (\alpha) = (\chi')^* \delta (\alpha)
= \delta (\alpha)({\mathbb T}^{Q_Y} )$$
Putting together the previous two displayed equalities, 
we get 
$$\pa_{(X,Y)}(\alpha(Q))= \delta (\alpha)({\mathbb T}^{Q_Y} )$$
as desired. This completes the proof of the equivalence of
Theorems \ref{main theorem} 
and \ref{universal gysin form of main theorem}. \hfill $\square$

\subsection{Standard models for degenerating quadrics}

Let $K'$ be the universal line bundle on $B\GG_m$. 
Consider any $n\ge 1$ and the universal triple $(E_{n-1},L_{n-1},b_{n-1})$ on 
the stack $BGO_{n-1}$. Consider the projections 
$p_1:BGO_{n-1}\times B\GG_m \to BGO_{n-1}$ and 
$p_2:BGO_{n-1}\times B\GG_m \to B\GG_m$.
Let $(\EE,\LL,\bb) = p_1^*(E_{n-1},L_{n-1},b_{n-1})$ and  
$K = p_2^*K'$ be the pullbacks.
Let the stack $\N$ be the total space of the line bundle $\LL$ on 
the stack $\Zz = BGO_{n-1}\times B\GG_m$, and let $\pi :\N \to \Zz$
be the projection $1$-morphism. 
We regard $\Zz$ as a closed substack $\Zz\hra \N$, embedded by 
the zero section. Note that the stacks $\Zz$ and $\N$ are smooth over $S$.

We now define a certain quadratic triple on $\N$ following [H-N-2],
called the {\bf model triple}. 
Let $\tau \in \Gamma(\N, \pi^*(\LL))$ be the tautological section,
which vanishes exactly on the zero section 
$BGO_{n-1}\times B\GG_m= \Zz \subset \N$,
with vanishing multiplicity $1$. Then we get a quadratic triple 
(where $\oplus$ denotes the orthogonal direct sum)
$$T = 
\left((\OO_{\N}, \pi^*\LL, \tau) \oplus \pi^*(\EE,\LL,\bb)\right)\otimes K$$
on $\N$ of dimension $n$, which has rank $n$ on $\N -\Zz$ and has
rank $n-1$ on $\Zz$, which means it is a mildly degenerating triple on 
$(\N, \Zz)$, which is nondegenerate on $\N-\Zz$ and 
minimally degenerate on the divisor $\Zz$, with degeneration 
multiplicity $\nu_{\Zz}(\det(T)) = 1$.

Let  $\varphi: \Zz \to \D$ be the classifying morphism 
of the degenerate triple 
$$T_0 = 
\left((\OO_{\Zz}, \LL, 0) \oplus (\EE,\LL,\bb)\right)\otimes K$$
which is the restriction of $T$ to $\Zz\hra \N$.
In the opposite direction, 
let the $1$-morphism $\psi: \D \to \Zz$ be defined as follows. 
For any $S$-scheme $Y$, an object of $\D(Y)$ is a minimally degenerate
quadratic triple $(F,L,q)$ on $Y$.  
Hence $\ker(q)$ is a line bundle on $Y$ so it is 
an object of $B\GG_m(Y)$, while
$(F/\ker(q), L, \ov{q})$ (where $\ov{q}$ is induced by $q$) 
is an object of $BGO_{n-1}(Y)$. 
As $\Zz = BGO_{n-1} \times B\GG_m$, we can define 
$$\psi(F,L,q) = ((F/\ker(q), L, \ov{q})\otimes \ker(q)^{-1},\, \ker(q))
\in Ob\,\Zz(Y).$$

\stm\label{rho equals p times psi}{\bf (Factorization of $\rho$) }
The $1$-morphism $\rho : \D \to BGO_{n-1}$, defined earlier by 
$\rho(F,L,q) = (F/\ker(q), L, \ov{q})\otimes \ker(q)^{-1}$, 
factors as
$$\rho = p_1\circ \psi$$
where $p_1 : BGO_{n-1} \times B\GG_m \to BGO_{n-1}$
is the projection.

The composite $1$-morphism $\varphi\circ\psi : \D \to \D$ is of importance 
to us. In terms of the above notation, we have 
$$\varphi\circ\psi (F,L,q) = 
(\ker(q), L, 0)\oplus (F/\ker(q), L, \ov{q}).$$

\begin{lemma}\label{A1 homotopy lemma}
The $1$-morphisms $\varphi: \Zz \to \D$ and $\psi: \D \to \Zz$
satisfy the following.

\noindent{(1)} $\psi\circ\varphi = \id_{\Zz}$.

\noindent{(2)} $\varphi\circ\psi : \D \to \D$ is 
$\AA^1$-homotopic to $\id_{\D}$, that is, there exists a $1$-morphism
$F : \D \times \AA^1 \to \D$ such that $\varphi\circ\psi = F_0$ 
and $\id_{\D} = F_1$, where $F_t$ denotes $F|_{\D\times \{ t\}}$ for 
$t\in \AA^1(\Z[1/2]) = \Z[1/2]$. Hence 
by Remark \ref{A1 homotopy invariance}, the cohomology map 
$\psi^*\circ \varphi^* : H^*(\D) \to H^*(\D)$ is identity.

\noindent{(3)} Consequently, the cohomology maps 
$\psi^* :  H^*(\Zz) \to H^*(\D)$ and $\varphi^* : H^*(\D) \to H^*(\Zz)$ 
are isomorphisms, which are inverse to each other.
\end{lemma}

\proof Let $0\to E'\to E \to E''\to 0$ is a short exact sequence of vector 
bundles on an $S$-scheme $Y$. Let $\AA^1$ be the affine line over 
with coordinate $t$. On $Y\times \AA^1$ 
consider the homomorphism $g: p_Y^*E' \to p_Y^*E'\oplus p_Y^*E$ 
which sends $v \mapsto (tv,v)$. Then $\F = \coker(g)$ fits
in a short exact sequence 
$0\to p_Y^*E' \to \F \to p_Y^*E'' \to 0$
of vector bundles on $Y\times \AA^1$, which specializes
to the given exact sequence on $Y$ for $t=1$ and specializes to the split
exact sequence $0\to E'\to E'\oplus E'' \to E''\to 0$ for $t=0$. 
The above construction is functorial, hence it also works 
for short exact sequences of vector bundles over algebraic stacks. 

Applying this to the short exact sequence 
$0\to \ker(b_n|_{\D}) \to \E_n|_{\D} \to \ov{\E_n|_{\D}} \to 0$ on 
$\D$ which arises from the universal minimally degenerate
triple $\T_{\D} = (\E_n,\L_n,\beta_n)|_{\D}$, 
we get a short exact sequence 
$0\to p_{\D}^*\ker(b_n|_{\D}) \to \F \to p_{\D}^*\ov{\E_n|_{\D}} \to 0$ 
of vector 
bundles on $\D \times \AA^1$. 
Let $\wh{b} : \F \otimes \F \to p_{\D}^*L$ be induced by 
$p_{\D}^*\ov{b}$. 
This defines a minimally degenerate triple 
$\wh{\T_{\D}} = (\F, p_{\D}^*L, \wh{b})$ on $\D \times \AA^1$.
Take $F  = \chi_{\wh{\T_{\D}}} : \D\times \AA^1 \to \D$
to be its classifying morphism. \qed

Let $\chi : \N \to \M$ be the classifying morphism of the 
mildly degenerating quadratic triple
$T = 
\left((\OO_{\N}, \pi^*\LL, \tau) \oplus \pi^*(\EE,\LL,\bb)\right)\otimes K$
on $(\N,\Zz)$ 
defined above. Under it, the pullback of $\D\subset \M$ is 
the zero section $\Zz \subset \N$, with multiplicity $1$. 
Both the pairs $(\M,\D)$ and $(\N,\Zz)$ admit
long exact Gysin sequences, and we have a commutative diagram 
as follows, where the vertical maps are induced by the classifying
morphism of $T$ and its restrictions.
$$\begin{array}{ccccccccc}
\ldots \to &H^i(\M)&\to&H^i(\M-\D)&\stackrel{\pa}{\to} & 
H^{i-1}(\D)&\to&H^{i+1}(\M)&\to\ldots\\
& ~~~ \da {\scriptstyle \chi^*}&  & ~~~ \da {\scriptstyle \chi|_{\N -\Zz}^*}& &
 ~~~ \da {\scriptstyle \chi|_{\Zz}^*}&& ~~~ \da {\scriptstyle \chi^*}& \\
\ldots \to & H^i(\N)&\to &H^i(\N-\Zz)&\stackrel{d}{\to}& 
H^{i-1}(\Zz)&\to&H^{i+1}(\N)&\to\ldots
\end{array}
$$
Note that $\M -\D$ is isomorphic to $BGO_n$, and  
$\N - \Zz$ is isomorphic to $BO_{n-1}\times B\GG_m$. Under
these isomorphisms, the above 
map $\chi|_{\N -\Zz}^* : H^i(\M-\D) \to H^i(\N -\Zz)$ is 
the cohomology map $H^i(BGO_n)\to H^i(BO_{n-1}\times B\GG_m)$
associated to the $1$-morphism 
$BO_{n-1}\times B\GG_m \to  BGO_n$ of stacks 
which is induced by the group scheme homomorphism
$$ %
V: O_{n-1}\times \GG_m \to GO_n : (g, \lambda) \mapsto 
\left(\begin{array}{cc} 1& \\ & g \end{array}\right)\lambda.$$
Hence the commutative diagram 
$$\begin{array}{ccc}
H^i(\M-\D)&\stackrel{\pa_{(\M,\D)}}{\to} & H^{i-1}(\D)\\
 ~~~ \da {\scriptstyle \chi|_{\N -\Zz}^*}&
& ~~~ \da {\scriptstyle \chi|_{\Zz}^*} \\
H^i(\N-\Zz)&\stackrel{\pa_{(\N,\Zz)}}{\to}& H^{i-1}(\Zz)\\
\end{array}
$$
becomes the commutative diagram
$$\begin{array}{ccc}
H^i(BGO_n)&\stackrel{\pa_{(\M,\D)}}{\to} &  H^{i-1}(\D)\\
{\scriptstyle B(V)^*} \da  ~~~& &  ~~~ \da {\scriptstyle \varphi^*} \\
H^i(BO_{n-1}\times B\GG_m) & \stackrel{d}{\to} &
                  H^{i-1}(BGO_{n-1} \times B\GG_m) \\
\end{array}
$$

\stm\label{factorization of gysin}{\bf (Factorization of the Gysin map
$\pa: H^i(\M-\D) \to  H^{i-1}(\D)$) } 
As $\varphi^*$ and $\psi^*$ are inverses by Lemma \ref{A1 homotopy lemma},  
the above commutative diagram gives the following commutative diagram.
$$
\begin{array}{ccc}
H^i(BGO_n)&\stackrel{\pa_{(\M,\D)}}{\to} &  H^{i-1}(\D)\\
{\scriptstyle B(V)^*} \da  ~~~& &  ~~~ \upa {\scriptstyle \psi^*} \\
H^i(BO_{n-1} \times B\GG_m) & \stackrel{d}{\to} &
                  H^{i-1}(BGO_{n-1} \times B\GG_m) \\
\end{array}
$$

\begin{lemma}\label{commutativity where primitivity matters}
Let $PH^*(BGO_n) \subset H^*(BGO_n)$ be the primitive subring. 
Then we have a commutative diagram
$$\begin{array}{ccc}
PH^i(BGO_n)&\hra &  H^i(BGO_n)\\
 ~~~ \da {\scriptstyle B(v)^*}& & ~~~ \da {\scriptstyle B(V)^*} \\
H^i(BO_{n-1}) & \stackrel{{p'_1}^*}{\to} & H^i(BO_{n-1}\times B\GG_m) \\
\end{array}
$$
where ${p'_1}^*$ is induced by the projection
$p'_1:BO_{n-1} \times B\GG_m \to BO_{n-1}$.
\end{lemma}

\proof Let $(F,q)$ be the universal quadratic form on $BO_{n-1}$, where
$F$ is a vector bundle of rank $n-1$ and $q$ is a nondegenerate quadratic
form on $F$ with values in $\OO_{BO_{n-1}}$. Let $K$ be the universal
line bundle on $B\GG_m$. 
The morphism $B(V):  BO_{n-1}\times B\GG_m \to BGO_n$
is the classifying map of the nondegenerate quadratic triple
$T = \left((\OO, \OO, 1)\oplus (F,\OO,q)\right)\otimes K$.
If $\alpha \in H^*(BGO_n)$ then 
$B(V)^*(\alpha) = \alpha(T)$. As $\alpha \in PH^*(BGO_n)$, we have
$\alpha(T) = \alpha (T \otimes K^{-1})
= \alpha ((\OO, \OO, 1)\oplus (F,\OO,q))
= B(v)^*(\alpha)$. 
\hfill$\square$

\subsection{Completion of the proof of the main theorem}

\stm\label{p and p' commutativity}
Let $p_1 : BGO_{n-1} \times  B\GG_m\to BGO_{n-1}$ and 
${p'}_1 : BO_{n-1}\times B\GG_m  \to BO_{n-1}$ be the projections.
Then the following diagram commutes, where the horizontal maps are
Gysin boundary maps.
$$\begin{array}{ccc}
H^i(BO_{n-1}) & \stackrel{d_{n-1}}{\to} & H^{i-1}(  BGO_{n-1})\\
{\scriptstyle {p'}_1^*} \da  ~~~& &~~~ \da {\scriptstyle p_1^*} \\
H^i(BO_{n-1}\times  B\GG_m) &\stackrel{d}{\to}&
                                   H^{i-1}(BGO_{n-1}\times  B\GG_m) 
\end{array}%
$$

As before, we identify $\M -\D$ with $BGO_n$, under which 
$PH^*(\M -\D)$ gets identified with $PH^*(BGO_n) \subset H^*(BGO_n)$.
We then have the following sequence of equalities.
\begin{eqnarray*}
\pa|_{PH^*(BGO_n)} 
& = & \psi^*\circ d\circ B(V)^* 
        \mbox{ by Statement \ref{factorization of gysin},}\\
& = & \psi^*\circ d \circ {p'}_1^* \circ B(v)^* 
        \mbox{ by Lemma \ref{commutativity where primitivity matters},}\\
& = & \psi^*\circ p_1^*\circ d_{n-1} \circ B(v)^*
   \mbox{ by Statement \ref{p and p' commutativity}},\\
& = & \rho^*\circ d_{n-1} \circ B(v)^*
   \mbox{ by Statement \ref{rho equals p times psi}}.
\end{eqnarray*}
Hence the Gysin boundary map 
$\pa_{(\M,\D)} : H^*(\M-\D)\to  H^{*-1}(\D)$ satisfies 
the formula 
$$\pa|_{PH^*(\M -\D)} = \rho^* \circ d_{n-1} \circ B(v)^*$$
which completes the proof of
Theorem \ref{universal gysin form of main theorem},
and hence of Theorem \ref{main theorem}. \qed

\subsection{The case of a separably closed base field $k = k_{sep}$}

When our base $S$ is of the form $\spec k$ for a
separably closed field $k$ of characteristic $\ne 2$, 
the rings $H^*(BGO_{n,k},\FF_2)$ and
$PH^*(BGO_{n,k},\FF_2)$ are the same as in the topological case by [Bh], so
the Theorem \ref{main theorem} now implies that
all the calculations in [H-N-2] apply without change over such a $k$.

In [Bh], it was shown that the key topological facts on which the 
the calculation of $H^*_{sing}(BGO_{2n}(\C),\FF_2)$ in [H-N-1]
is based have suitable analogs in 
the algebraic category over a separably
closed base field $k$ of characteristic $\ne 2$, and hence we get the 
following determination of $H^*(BGO_{n,k},\FF_2)$.

First, recall that the cohomology ring $H^*(BGL_{n,k},\FF_2)$ 
is the polynomial ring 
$\FF_2[\ov{c}_1,\ldots, \ov{c}_n]$ in $n$ variables
$\ov{c}_i\in H^{2i}(BGL_{n,k},\FF_2)$ which are the 
{\bf universal mod $2$ Chern classes}. Also, 
the cohomology ring $H^*(BO_{n,k},\FF_2)$ 
is the polynomial ring
$\FF_2[w_1,\ldots, w_n]$ in $n$ variables
$w_i\in H^i(BGL_{n,k},\FF_2)$ which 
are the {\bf universal Stiefel-Whitney classes} 
of Laborde [La] over $k=k_{sep}$ (see, for example, [E-K-V]). 
Under the homomorphism
$H^*(BGL_{n,k},\FF_2)\to H^*(BO_{n,k},\FF_2)$ induced by the inclusion
$O_{n,k}\hra GL_{n,k}$, we have $\ov{c}_i\mapsto w_i^2$.

\stm{\bf The rings $H^*(BGO_{2n+1,k},\FF_2)$.}
For the odd special orthogonal group scheme $SO_{2n+1,k}$
we have $H^*(BSO_{2n+1,k},\FF_2) = \FF_2[w_2,\ldots, w_{2n+1}]$, 
as follows from the direct product decomposition 
$O_{2n+1} = \mu_2 \times SO_{2n+1,k}$ by using the K\"unneth formula
(see section 5 of [Bh] for details of the argument).
We have an isomorphism 
$\GG_m\times SO_{2n+1} \to GO_{2n+1}$ defined in terms of valued points by
$(\lambda, g)\mapsto \lambda g$. It follows that
$BGO_{2n+1} = B\GG_m\times BSO_{2n+1}$. Hence we have  
$$H^*(BGO_{2n+1,k},\FF_2) = \FF_2[c, w_2,\ldots, w_{2n+1}]$$
where $c$ denotes the pullback of the first Chern class
$\ov{c}_1$ from $H^2(BG_{m,k},\FF_2)$, 
and the $w_2,\ldots, w_{2n+1}$ are the pull-backs 
of the Hasse-Witt classes in 
$H^*(BSO_{2n+1},\FF_2)$. Note that $c, w_2,\ldots, w_{2n+1}$ 
are algebraically independent over $\FF_2$,
and are homogeneous of degrees $2, 2, 3, \ldots, 2n+1$.

\stm{\bf The rings $H^*(BGO_{2n,k},\FF_2)$.} 
This ring is generated by the cohomology classes 
$\lambda$, $a_{2i-1}$, $b_{4j}$, and $d_T$ over $\FF_2$, 
where $i, j \in \{ 1, \ldots, n\}$ and $T\subset \{ 1, \ldots, n\}$
has cardinality $|T| \ge 2$. These are defined as follows. 
Let $(E_{2n},L_{2n},b_{2n})$ denote the universal quadratic triple on $BGO_{2n}$.
Then $\lambda = \ov{c}_1(L_{2n})$ and $b_{4j}= \ov{c}_{2j}(E_{2n})$.
Next, consider the Gysin boundary map $d : H^*(BO_{n,k},\FF_2) \to
H^{*-1}(BGO_{n,k},\FF_2)$ that we introduced earlier.
Then 
$a_{2i-1} = d(w_{2i}) \in H^{2i-1}(BGO_{n,k},\FF_2)$, 
and for any
$T = \{ i_1 ,\ldots, i_r\}\subset \{ 1 ,\ldots, n \}$, we have
$d_T = d(w_{2i_1}\cdots w_{2i_r}) \in 
H^{2(i_1 + \ldots + i_r) -1}(BGO_{n,k},\FF_2)$. 
These generators satisfy relations given explicitly in [H-N-1] Theorem 3.9.

\stm {\bf Relation with Stiefel-Whitney and Chern classes.}
When $k = k_{sep}$, 
the ring homomorphisms 
$H^*(BGL_{n,k},\FF_2) \to H^*(BGO_{n,k},\FF_2) \to 
H^*(BO_{n,k},\FF_2)$,
induced by the inclusions $O_n\hra GO_n\hra GL_n$,
are given by the same explicit 
formulas in terms of the generators of $H^*(BGO_{n,k},\FF_2)$
which are deduced in Section 3 of [H-N-2] in the 
topological category for $H^*_{sing}(BGO(\C),\FF_2)$.

\stm 
With the above explicit descriptions of $H^*(BGO_{n,k},\FF_2)$ in the even and
odd cases, all the calculations Section 8 of [H-N-2] hold in the 
algebraic category over such a $k$. 
These include Corollaries 8.1 and 8.2 which give formulas 
in terms of generators and relations for the 
`even rank degenerating to odd rank' and `odd rank degenerating to even rank'
cases of the Theorem  \ref{main theorem}, illustrated with complete lists of
generators of the primitive rings $PH^*(BGO_{n,k},\FF_2)$
and their images under $\pa$ in all ranks up to $n=6$.

\rem\label{when kummer simplifies}
For any scheme $X$ over $\Z[1/2]$ and any 
$y \in \Gamma(X,\OO_X)^{\times}$, let $(y) \in H^1(X,\FF_2)$
denote the Kummer class of $y$. If $X$ lies over a 
separably closed field base field of characteristic $\ne 2$,
then note that $(y)^2 = 0\in H^2(X,\FF_2)$, and moreover, 
$(-1) = 0 \in H^2(X,\FF_2)$. These properties of the Kummer class also hold
in the topological category for Kummer classes of 
nowhere vanishing complex valued continuous functions $y$. 
These properties, which are crucially used in [H-N-1] and [H-N-2],
are not available over an arbitrary base such as $\Z[1/2]$. In fact 
all powers $\kappa^n$ of the Kummer class 
$\kappa = (-1) \in H^1(\Z[1/2],\FF_2)$  
are nonzero, as may be seen by base change to $\R$. 
This is one of the reasons why the explicit determination of the rings
$H^*(BGO_n,\FF_2)$ and $PH^*(BGO_n,\FF_2)$ in terms of generators and relations
is more complicated over a general base.
This question will be addressed elsewhere.

\parindent=0pt 
\parskip=4pt


{\small

\section*{References} \addcontentsline{toc}{section}{References}

[Be] Behrend, K.A. : Derived $\ell$-adic categories for algebraic stacks. 
Mem. Amer. Math. Soc. 163 (2003), no. 774

[Bh] Bhaumik, S. : Characteristic classes for $GO(2n)$ in \'etale
cohomology. Preprint, 2012. To appear in Proc. Indian Acad. Sci. 
(http://arxiv.org/abs/1201.4628) 

[EGA-IV$_2$] Grothendieck, A. : \'El\'ements de g\'eom\'etrie alg\'ebrique. IV. 
\'Etude locale des sch\'emas et des morphismes de sch\'emas, seconde partie. 
Inst. Hautes \'Etudes Sci. Publ. Math. No. 24, 1964.

[E-K-V] Esnault, H., Kahn, B., Viehweg, E. : 
Coverings with odd ramification and Stiefel-Whitney classes. 
J. Reine Angew. Math. 441 (1993), 145-188.

[Fu] Fujiwara, K. :
A proof of the absolute purity conjecture (after Gabber). 
Algebraic geometry 2000, Azumino (Hotaka), 153-183, 
Adv. Stud. Pure Math., 36, Math. Soc. Japan, 2002. 




[Ju] Jouanolou, J. P. : 
Une suite exacte de Mayer-Vietoris en K-th\'eorie alg\'ebrique. 
Algebraic K-theory, I: Higher K-theories 
(Proc. Conf., Battelle Memorial Inst., Seattle, 1972), 
293-316. Lecture Notes in Math., 341, Springer, 1973.

[H-N-1] Holla, Y.I. and Nitsure, N. : 
Characteristic classes for $GO(2n, \C)$, Asian Journal
of Mathematics 5, 2001, 169-182.

[H-N-2] Holla, Y.I. and Nitsure, N. : 
Topology of quadric bundles. Internat. J. Math. 12 (2001), 1005-1047.

[La] Laborde, O. :
Classes de Stiefel-Whitney en cohomologie étale. 
Colloque sur les Formes Quadratiques (Montpellier, 1975).
Bull. Soc. Math. France Suppl. Mem. No. 48 (1976), 47–51. 

[L-MB] Laumon, G. and Moret-Bailly, L. : 
{\it Champs algébriques}. Springer 2000.

[L-O] Laszlo, Y. and Olsson, M. : The six operations for sheaves on 
Artin stacks I: Finite coefficients.  
Publ. Math. IH\'ES 107 (2008), 109-168. 

[Na-Ra] Narasimhan, M.S. and Ramanan, S. : Geometry of Hecke cycles -I. 
{\it C. P. Ramanujam - a tribute}, 291-345, 
TIFR Studies in Mathematics 8, Springer, 1978.

[N-1] Nitsure, N. : Topology of conic bundles. 
J. London Math. Soc. (2) 35 (1987) 18-28. 
 
[N-2] Nitsure, N. : Cohomology of desingularization of moduli space
of vector bundles. Compositio Math. 69 (1989) 303-339.

[O] Olsson, M. : Sheaves on Artin stacks. 
J. Reine Angew. Math. 603 (2007), 55-112.

[Se] Seshadri, C.S. : Geometric reductivity over arbitrary base. 
Advances in Math. 26 (1977), 225-274. 

[SGA 4] Artin, M., Grothendieck, A., Verdier, J.-L. (ed.) : 
{\it Th\'eorie des topos et cohomologie \'etale des sch\'mas}, 
Sem. Geom. Alg. Bois Marie 1963-64, 
LNM 269, 270 and 305, Springer 1972-73.

[SGA $4{1\over 2}$]  Deligne, P. : {\it Cohomologie \'etale}.
(With the collaboration of J. F. Boutot, \\
A. Grothendieck, L. Illusie and J.-L. Verdier.) LNM 569, Springer 1977.

[Td] Toda, H. : Cohomology of Classifying Spaces. In the volume
{\sl Homotopy Theory and Related Topics} (Editor H. Toda),   
Advances in Pure Math. {\bf 9} (1986) 75-108.

[Tot] Totaro, B. : Chow ring of a classifying space, 
{\it Algebraic K-theory} (Seattle 1997) 249-281, 
Proc. Symp. in Pure Math. 67,
Amer. Math. Soc., 1999.

\bigskip 




Saurav Bhaumik $<$saurav@math.tifr.res.in$>$ and
Nitin Nitsure $<$nitsure@math.tifr.res.in$>$

Tata Institute of Fundamental Research,
Homi Bhabha Road, Mumbai 400 005, India.

\bigskip 

\centerline{01 February 2013, revised on 24 April 2013}

}

\end{document}